%
% CCH paper "The 2x2 Matrix Mortality Problem with at Most One Invertible Matrix is Decidable"
%
% Run twice to get cross-references correct
%

\def \introduc{Section~1}
\def \mainres{Section~2}
\def \one{Theorem~1}
\def \product{Lemma~2}
\def \cosine{Lemma~3}
\def \power{Lemma~4}
\def \mmInv{Section~3}
\def \two{Theorem~5}
\def \refs{Section~4}

\def \bb{\hbox {\rm [1]}}
\def \twentyeighteen{\hbox {\rm [2]}}
\def \chthm{\hbox {\rm [3]}}
\def \nucrod{\hbox {\rm [4]}}
\def \pater{\hbox {\rm [5]}}
\def \shank{\hbox {\rm [6]}}
\def \twentyfive{\hbox {\rm [7]}}

\def\BandB{Bournez and Branicky}
\def\rank{\mathop{\rm rank}}
\def\adj{\mathop{\rm adj}}
\long\def\ignore#1{}

\let\UseNewFonts N
%
% Euler font version
%

\newdimen\fudgeAmount
\def\fudge#1!{\fudgeAmount #1
	\advance\hsize by 2\fudgeAmount \advance\vsize by 2\fudgeAmount 
	\advance\hoffset by -\fudgeAmount \advance\voffset by -\fudgeAmount
	}

%\font\textfont=cmr8 

%%% stuff to use Euler fonts as default math fonts
\font\teneurm=eurm10
\font\seveneurm=eurm7
\font\fiveeurm=eurm5
\newfam\eurmfam
\textfont\eurmfam=\teneurm
\scriptfont\eurmfam=\fiveeurm
\scriptscriptfont\eurmfam=\fiveeurm

% automatic Euler for math mode
\textfont1=\teneurm
\scriptfont1=\seveneurm
\scriptscriptfont1=\fiveeurm
% new defs, such as \vec
\newdimen\vecArgWidth
\newdimen\vecArgHeight
\newdimen\vecEps
\newbox\vecArgBox
\def\Rightarrowfill{$\scriptscriptstyle\mathsurround0pt 
	\smash- \mkern-7mu\cleaders\hbox{$\mkern-2mu 
	\smash- \mkern-2mu$}\hfill\mkern-7mu \mathord\rightarrow$}
\def\vec#1{
	\setbox\vecArgBox\hbox{$#1$}
	\vecArgWidth \wd\vecArgBox 
	\vecArgHeight \ht\vecArgBox
	\vecEps\vecArgWidth \divide\vecEps by 50
	\advance\vecArgHeight by 0.5pt
	\rlap{\raise\vecArgHeight\hbox to\vecArgWidth{\Rightarrowfill}}
	\hbox{$#1$}
	}

\ifx\UseNewFonts Y
    % \UseNewFonts; hyperbold, etc.
    
    \font\bbold bbold10

\fi

\let\PAR\par
\ifx\ShowSolns Y

	\def\ans#1{\setbox\ansbox\hbox{\rm #1}\ansdepth\dp\ansbox\advance\ansdepth by 3.5pt{\lower\ansdepth \fbox{\box\ansbox}}}
	\def\mathans#1{\setbox\ansbox\hbox{$#1$}\ansdepth\dp\ansbox\advance\ansdepth by 3.5pt{\lower\ansdepth \fbox{\box\ansbox}}}
	\def\ansSameLine{\hfill \quicksoln}
	\def\pp#1{\hbox{$+#1$~points}}
	\def\np#1{\hbox{$-#1$~points}}

	\long\def\soln#1{\bigskip{\narrower\noindent\ital{Solution:}\sl\ #1 \bigskip\vfill}}
	\long\def\grading#1:#2\endgrading{\bigskip\centerline{\hfill 
		$\left[ ~~\vcenter{\noindent\ital{#1}:\sl\ #2 \PAR}~~\right]$\hfill}\bigskip}
	\def\hint#1{\bigskip{\narrower\noindent\ital{Hint:}\sl\ #1 \bigskip\vfill}}
	\def\quicksoln#1{\ital{Solution:}{\sl\ #1}}
	\def\quickans#1{\medskip\itemitem{}\ital{Answer:} #1 \vfill}
	\def\quickmathans#1{\mathans{#1} \vfill}
	\def\veryQuickAnsRlap#1{\rlap{\veryQuickAns{#1}}}
	\def\veryQuickAns#1{{\sl #1}}
	\def\ignore{}
	
		% unless changed later
\else

	\def\ans#1{}
	\def\mathans#1{}
	\def\ansSameLine#1{}
	\def\pp#1{}
	\def\np#1{}

	\long\def\soln#1{\bigskip\vfil}
	\long\def\grading#1\endgrading{}
	\def\hint#1{\bigskip}
	\def\quicksoln#1{}
	\def\quickans#1{}
	\def\quickmathans#1{}
	\def\veryQuickAns#1{}
	\def\veryQuickAnsRlap#1{}
	\long\def\ignore#1{}
	
\fi

\def\dst{\displaystyle}

\def\resize#1{\magnification=#1
\expandafter\ifx\csname pdfpageheight\endcsname\relax \else \pdfpageheight 11true in \pdfpagewidth 8.5true in
\advance\hoffset by -0.2true in \advance\voffset by -0.2true in
\fi
}

\long\def\ignore#1{}
\everypar{\looseness-1 }
\everymath{\displaystyle}

\newcount\itemmno \itemmno0
\newcount\itemmmno \itemmmno0
\newcount\itemno \itemno0
\def\itemitemitem{\par\indent\indent \hangindent3\parindent \textindent}

\def\iitem{\global\itemmno0 \global\advance\itemno by 1  \medskip\item{\the\itemno.} \global\itemmno 0 \global\itemmmno 0}
\def\iitemm{\global\advance \itemmno by 1  \itemitem{{\rm\ifcase\itemmno\or a\or b\or c\or d\or e\or f\or g\or h\or i\or j\or k\or$\ell$\or m\or n\fi.}}\global\itemmmno0 }
\def\iitemmm{\global\advance \itemmmno by 1  \itemitemitem{\ifcase\itemmmno\or i\or ii\or iii\or iv\or v\or vi\or vii\or viii\or ix\fi.} }
\def\bonus{\itemmno0 \global\advance\itemno by 1  \medskip\item{*~\the\itemno.} \global\itemmno 0 (Optional) }
\def\bi{\item{$\bullet$}\strut }
\def\bii{\itemitem{$\bullet$}\strut }

\def\MCfour#1\\#2#3#4#5{\iitem #1\hfill\rlap{\vrule width0.5in depth0.5pt height0pt}\medskip
\itemitem{A.} #2 \smallskip 
\itemitem{B.} #3 \smallskip 
\itemitem{C.} #4 \smallskip 
\itemitem{D.} #5 \smallskip 
\vfill}

\def\MCfive#1\\#2#3#4#5#6{\iitem #1\hfill\rlap{\vrule width0.5in depth0.5pt height0pt}\medskip
\itemitem{A.} #2 \smallskip 
\itemitem{B.} #3 \smallskip 
\itemitem{C.} #4 \smallskip 
\itemitem{D.} #5 \smallskip 
\itemitem{E.} #6 \smallskip
\vfill}

\def\LT(#1){{\cal L}\(#1\)} % Laplace transform
\long\def\fbox#1{\vbox{\hrule\hbox{\vrule\kern3pt\vbox{\kern3pt#1\kern3pt}\kern3pt\vrule}\hrule}}
\def\twoonline#1#2{\par\line{\indent \hbox to 2.5in{#1 \hfill} \hbox to 2.5in{#2
	\hfill}\hfill}}

\newbox\ansbox
\newdimen\ansdepth

%%%%%%%%%%%%%%%%%%% MAT 27X:

\def\ital#1{{\it #1\/}}

\def\set#1{{\left\{ #1 \right\} }}
\def\size#1{{\left| #1 \right| }}

\def\roots#1\of#2{\root\scriptstyle #1\of{#2}}

\def\(#1\){\left(\dst #1\right)}
\def\[#1\]{\left[\dst #1\right]}

\def\R{\hbox{\bbold R}}
\def\N{\hbox{\bbold N}}

\def\calP(#1){{\cal P}(#1)}  %%% power set

\def\pmb#1{\setbox0=\hbox{#1}%
\kern-.025em\copy0\kern-\wd0
\kern.05em\copy0\kern-\wd0
\kern-.025em\raise.0233em\box0 }

\def\doubleint_#1{\mathop{\int\!\!\!\int}_{\dst #1 ~~}}
\def\tripleint_#1{\mathop{\int\!\!\!\int\!\!\!\int}\limits_{\dst #1~~}}
\def\part#1#2{{\partial #1\over \partial #2}}

\def\adots{\mathinner{\mkern2mu\raise1pt\hbox{.}\mkern2mu\raise4pt\hbox{.}%
	\mkern2mu\raise7pt\hbox{.}\mkern1mu}}

%%%%%%%%%%%%%%%%%%% MAT 275:

\def\ltsim{\rlap{\lower10pt\hbox{~$\widetilde{~~}$}}<}
\def\gtsim{\rlap{\lower10pt\hbox{~$\widetilde{~~}$}}>}

%%%%%%%%%%%%%%%%%%% MAT 242:

\def\geogebra{\hbox{\rlap{$\bigcirc$}\hskip2pt\raise.5pt\hbox{\hsize=\bigcircwd \hfill $\scriptstyle\rm G$ \hfill}}}  
\def\row#1{\rlap{$\bigcirc$}\hskip3pt\raise.5pt\hbox{\hsize=\bigcircwd \hfill $\scriptstyle #1$ \hfill}}  
		% \row2 ; \row3\, + 5 \row1

\newdimen\bigcircwd \setbox1\hbox{$\bigcirc$} \bigcircwd\wd1
\def\dst{\displaystyle}

% \matrix, with entries aligned flush right.
\def\rmatrix#1{\null\,\vcenter{\normalbaselines\mathsurround=0pt 
	\ialign{\hfil$##$&&\quad \hfil$##$\crcr
	\mathstrut\crcr\noalign{\kern-\baselineskip}
	#1\crcr\mathstrut\crcr\noalign{\kern-\baselineskip}}}\,}

% \matrix, without widely spaced columns.
\def\mmatrix#1{\null\,\vcenter{\normalbaselines\mathsurround=0pt 
	\ialign{\hfil$##$\hfil&&\hfil$##$\hfil\crcr
	\mathstrut\crcr\noalign{\kern-\baselineskip}
	#1\crcr\mathstrut\crcr\noalign{\kern-\baselineskip}}}\,}

% \matrix, without widely spaced columns, flush right.
\def\rmmatrix#1{\null\,\vcenter{\normalbaselines\mathsurround=0pt 
	\ialign{\hfil$##$&&\hskip0.0625in\hfil$##$\crcr
	\mathstrut\crcr\noalign{\kern-\baselineskip}
	#1\crcr\mathstrut\crcr\noalign{\kern-\baselineskip}}}\,}

% \matrix, without widely spaced columns, flush left.
\def\lmatrix#1{\null\,\vcenter{\normalbaselines\mathsurround=0pt 
	\ialign{$##$\hfil&&\hskip0.0625in$##$\hfil\crcr
	\mathstrut\crcr\noalign{\kern-\baselineskip}
	#1\crcr\mathstrut\crcr\noalign{\kern-\baselineskip}}}\,}

\def\brcmat#1{\left[\matrix{#1}\right]}

%%%%%%%%%%%%%%%%%%%% synthetic division macros

%%%%%%%%%%%%%%%%%%%% long division macros

 % phantom +

% Divide the polynomial ${x^3+2x^2-5}$ by ${3x+1}$.
%
%$$\longdiv{&&{1\over3}x^2 &+& {5\over9}x &-& {5\over27}\cr
%	\noalign{\vskip-8pt}&\multispan8\hrulefill\cr\noalign{\vskip-4pt}
%	3x+1 &\vrule height12pt& x^3 &+& 2x^2 &&& - &5\cr
%	\noalign{\vskip-11pt}\multispan2\hrulefill\cr
%	&& x^3 &+& {1\over3}x^2\cr
%	\noalign{\vskip-8pt}&&\multispan3\hrulefill\cr
%	&&&&{5\over3}x^2\cr
%	&&&&{5\over3}x^2 &+& {5\over9}x\cr
%	\noalign{\vskip-8pt}&&&&\multispan3\hrulefill\cr\noalign{\vskip-1pt}
%	&&&&&-&{5\over9}x &-& 5\cr 
%	&&&&&-&{5\over9}x &-& {5\over27}\cr
%	\noalign{\vskip-8pt}&&&&&&\multispan3\hrulefill\cr
%	&&&&&&&-&{130\over27}\cr}$$

%%%%%%%%%%%%%%%%%%%%%%%
% \eqalignnno: LHS, RHS, MIDDLE, LHS, RHS, RIGHT
% $$\eqalignnno{x &= y & so & y &= x & (eq 1)\cr  z &= w\cr  x+1 &= x+1 && x &= x \cr}$$
%
%%%%%%%%%%%%%%%%%%%%%%%

\newdimen\jot \jot=3pt

\catcode`\@=11
\newdimen\dimen@
\def\eqalignnno#1{\displ@y \tabskip=\centering
  \halign to \displaywidth{\hfil$\@lign\displaystyle{##}$\tabskip=0pt
    &$\@lign\displaystyle{{}##}$\hfil\tabskip=\centering
    &\hfil ##\hfil&\hfil$\@lign\displaystyle{##}$\tabskip=0pt
    &$\@lign\displaystyle{{}##}$\hfil\tabskip=\centering
    &\llap{$\@lign##$}\tabskip=0pt\crcr
    #1\crcr}}
\def\@lign{\tabskip=0pt\everycr={}} % restore inside \displ@y
\def\displ@y{\global\dt@ptrue \openup1\jot \m@th
  \everycr{\noalign{\ifdt@p \global\dt@pfalse \ifdim\prevdepth>-1000pt
      \vskip-\lineskiplimit \vskip\normallineskiplimit \fi
      \else \penalty\interdisplaylinepenalty \fi}}}
\newif\ifdt@p
\def\openup{\afterassignment\@penup\dimen@=}
\def\@penup{\advance\lineskip\dimen@
  \advance\baselineskip\dimen@ \advance\lineskiplimit\dimen@}
\def\m@th{\mathsurround=0pt }

%%%%%%%%%%%%%%%%%%%%%%%%%%%%%%%%
% Therefore macro, and \ergo (which is short for \cr \thbar\cr \therefore & )

\def\therefore{\hbox{.\hskip-1pt\raise.75ex\hbox{.}\hskip-1pt.}~}

\newcount\stepnum

\def\autoproofcheap#1{\global\stepnum0\halign{\global\advance\stepnum by 1\hbox to 0.5in{\hfil\strut(\the\stepnum)}~~$##$\quad\hfil&##\strut\hfil\cr #1}}
\def\autoproofcheapcheap#1{\global\stepnum0\halign{\global\advance\stepnum by 1\hbox to 0.285in{\hfil\strut(\the\stepnum)}~~$##$\quad\hfil&##\strut\hfil\cr #1}}

\font\bbold bbold10

\everypar={\looseness=-1}
\def\lastmod{2019 September 6}

\long\def\abstract#1{\medskip{\noindent\narrower Abstract: \quad #1\par} \medskip}

\newcount\secnum
\newcount\result
\newcount\fignum
\secnum0 
\fignum0
\result0
\def\newsection#1#2\par{\allowbreak{\global\advance\secnum by 1 
	\immediate\write1{\def\string#1{Section\noexpand~\the\secnum}}
	\bigskip\titlefont
	%\textfont1=\titlemath
	\noindent\llap{\tinybold\if\tokennames1\string#1'' 
	\fi}\hbox{}\the\secnum. #2
	\par}
	\if\tokennames1\headline{\hfill\rm\ifnum\pageno>1 #2 \fi\hfill}\fi
	\nobreak\rm\nobreak\medskip\noindent}
\long\def\newthm#1#2{\global\advance\result by 1 
	\immediate\write1{\def\string#1{Theorem\noexpand~\the\result}}
	\medskip\noindent\llap{\if\tokennames1\tinybold\string#1'' \fi}\bf Theorem 
	\the\result.\ \sl #2\rm\medskip}
\long\def\newlemma#1#2{\global\advance\result by 1 
	\immediate\write1{\def\string#1{Lemma\noexpand~\the\result}}
	\medskip
	%\vbox
	{\noindent\llap{\tinybold\if\tokennames1\string#1'' \fi}\bf Lemma 
	\the\result.\ \sl #2\medskip}}
\long\def\newconj#1#2{\global\advance\result by 1 
	\immediate\write1{\def\string#1{Conjecture\noexpand~\the\result}}
	\medskip\noindent\llap{\if\tokennames1\tinybold\string#1'' \fi}\bf Conjecture 
	\the\result.\ \sl #2\rm\par}

\font\titlefont cmr10 at 14pt
\font\bigtt=cmtt10 at 15pt
\font\tinybold=cmb10 at 7 pt	% cmb7
\font\titlemath cmsy10 at 14pt

\def\tokennames{0}	% 0 to omit marginalia & headers

\openin2=mortality.aux
\ifeof2\else\input mortality.aux \fi
\closein2

\openin2=mortalityrefs.aux
\ifeof2\else\input mortalityrefs.aux \fi
\closein2

\immediate\openout1=mortality.aux

\tokennames1

\def\dolist{\afterassignment\dodolist\let\next= }
\def\dodolist{\ifx\next\endlist \let\next\relax \else \\\let\next\dolist \fi \next}
\def\endlist{\endlist}
\def\\{\expandafter\ifcase\mode\ifx\next >\mode1\rangle\else\ifx \next< \mode1\langle
		\else\ifx\next[ \mode1 [\else {^\next}\fi\fi\fi
	\or\ifx\next] ] \mode2\else\ifx\next> \rangle \mode2\else\ifx\next< \langle \mode2
		\else \next\fi\fi\fi
	\or{^\next}\fi}

%%%%%%%%%%%%%%%%%%%%%%%%%%%%%%%%%%%%%%%%%%%%%%%%%%%%%%%%%%%%%%%%%%%%%%%%%%%%%%%%%%%%%%%%%%%%%%%%%%%%%%%%%%%%%%%%%5
%%%%%%%%%%%%%%%%%%%%%%%%%%%%%%%%%%%%%%%%%%%%%%%%%%%%%%%%%%%%%%%%%%%%%%%%%%%%%%%%%%%%%%%%%%%%%%%%%%%%%%%%%%%%%%%%%5
	% now the actual paper.
%\phantom{this}
%\vfil

\centerline{\titlefont The 2$\times\textfont2=\titlemath$2 Matrix Mortality Problem and Invertible Matrices}

\bigskip

\centerline{\titlefont Christopher Carl Heckman}\smallskip
\centerline{\bigtt Christopher.Heckman@asu.edu}\bigskip

\centerline{\titlefont School of Mathematical and Statistical Sciences}\smallskip
\centerline{\titlefont Arizona State University, Tempe, AZ, 85287--1804}\bigskip

\abstract{By modifying the proof of a paper by O.~Bournez and M.~Branicky, we establish that the Matrix Mortality
Problem is decidable with any finite set of $2\times2$ matrices which has at most one invertible matrix.
The same modification also shows that the number of non-invertible matrices is irrelevant.
}

\newsection\introduc{Introduction}

A set $S$ of $m$ matrices, each with dimensions $n\times n$, 
is mortal if there is some product of its entries (possibly with repetition)
is the zero matrix. The Matrix Mortality Problem [MMP] is to determine under which circumstances this problem is decidable,
whether there exists an algorithm that solves it. Here, all matrices will be assumed to be rational.

A previous result \pater\ establishes that the MMP is undecidable if
$n=3$; undecidability was decided for $(m,n) \ge (3,6),(5,4),(9,3),(15,2)$ by \twentyeighteen.

At the other extreme, the MMP is trivially decidable if $n=1$, or if $m=1$, or if the matrices are all upper-triangular. 
Bournez and Branicky also showed that the MMP is decidable if $m=n=2$. 
Decidability of the MMP was extended in \nucrod\ 
to sets of \ital{integral} matrices, each of whose determinants is $0$ or $\pm1$.

Henceforth, it will be assumed that $n=2$; all matrices involved will be $2\times2$.

\newsection\mainres{Main Result}

The decidability of the MMP with \hbox{$m=2$} is Theorem~4 of \bb. 
\BandB\ provide two proofs of their Theorem~4; they discovered a longer on on their own, and a referee discovered
a short but not-self-contained proof. It will be shown how to alter both proofs to obtain a proof of 
\one\ (below), even though one would be sufficient.

\newthm\one{The MMP is decidable if the set contains at most one invertible matrix.}

To prove this extension, some results are required. The first proves a bit more about what desired products will look
like.

\newlemma\product{{\bf(Lemma~2, \bb).} A finite set $F=\set{A_1,\ldots,A_m}$ of $2\times2$ matrices is mortal if and only if there exist
an integer $k$ and integers $i_1,\ldots,i_k\in\set{1,\ldots,m}$ with $A_{i_1}\cdots A_{i_k}=0$, and
\smallskip
\itemitem{\rm1.} $\rank A_{i_j} = 2$, for $1 < j < k$,
\itemitem{\rm2.} $\rank A_{i_j} < 2$, for $j\in \set{1,k}$.}

% this is also true for any number of matrices, as long as the noninvertible ones all have rank 1.

If the set $S$ of matrices has no invertible matrix, this will imply that $S$ is mortal iff 
$\set{AB \mid A,B\in S}$ contains the zero matrix, which can be checked in finite time. Otherwise, we will
write \hbox{$S=\set{A_1,\ldots,A_{m-1},B}$} where $B$ is the only invertible matrix in $S$. We also eliminate the possibility
that $A_i$ is the zero matrix for some $i$, making $\rank A_i=1$ for all $i$, and the possibility that $m=1$, in
which case $S$ is not mortal.

This amounts to checking whether $A_i B^k A_j=0$ for some nonnegative integers $i,j,k$. The case $i=j$ is handled
in \bb; the objective is to show that this equation can be effectively solved if $i\not=j$. Without loss of generality,
we assume that $i=1$ and $j=2$.

\medskip

\ital{Proof \#1 of \one.} The referee's proof. To every word $w = w_1 w_2\cdots w_n \in\set{B,A_1,A_2}^*$ we will 
associate the matrix
$A_w=w_1 w_2 \cdots w_n$. Let the language $Z$ consist of all words $w\in \set{B,A_1,A_2}^*$ such that $A_w=0$.
Then \one\ follows if one can effectively test if $Z$ is empty. In fact, the stronger claim that
$Z$ can be effectively computed will be proven.

Now write $A_1=ac^\top$ and $A_2=bd^\top$. By \product, if there exists a mortal product, then
$A_1 B^k A_2=0$ for some integer $k$. This condition is equivalent to $c^\top B^k b = 0$. The referee's proof
continues the same from this point; this ends Proof~\#1.

\medskip

\ital{Proof \#2 of \one.} The co-authors' proof. Once again, the problem reduces to determining whether there is an integer
$n$ such that $A_1 B^n A_2=0$. We check this relation algebraically using the Jordan forms of the matrices $A_1,A_2,B$.
Write:
$$A_1 = P_1^{-1} J_1 P_1, \quad B = P_2^{-1} J_2 P_2, \quad A_2 = P_3^{-1} J_3 P_3.$$ 
(We are trying to mimic the notation of the original paper as much as possible.)
Because $A_1$ and $A_2$ each have rank~1, we have
$$J_1 = \brcmat{\kappa_1 &0\cr 0&0\cr}, \quad J_3 = \brcmat{\kappa_2 &0\cr 0&0\cr},$$
for some nonzero number $\kappa_1$ (resp.~$\kappa_2$), which is rational because $A_1$ (resp.~$A_2$) 
is a rational matrix whose trace
equals $\kappa_1$ (resp.~$\kappa_2$).
$B$ is invertible, which makes 
$$J_2 = \brcmat{\lambda&0\cr 0&\mu\cr} \quad {\rm or} \quad \brcmat{\lambda&1\cr 0&\lambda\cr}.$$

Substituting the decompositions into the equation $A_1 B^n A_2=0$ yields the equation
$$P_1^{-1} J_1^{~} P_1^{~} P_2^{-1} J_2^n P_2^{~} P_3^{-1} J_3 P_3 =0,$$
which is equivalent to
$$J_1^{~} P J_2^n Q J_3=0,$$
since $P_1$ and $P_3$ are invertible, and
once we write $P=P_1^{~} P_2^{-1} = \brcmat{p&q\cr r'&p'\cr}$ and $Q=P_2 P_3^{-1} = \brcmat{s&q'\cr -r&s'\cr}$.
(The reason for choosing the entries this way will become clear soon.)

Once we make these substitutions, we find that the problem is now equivalent to testing whether
there is an integer $n$ such that
\medskip
\bii $\lambda^n p s - \mu^n q r=0$, when $J_2$ is of the first form; or
\smallskip
\bii $(ps-qr)\lambda-n p r = 0$, when $J_2$ is of the second form.
\medskip

(Note that there is a mistake in \bb's proof of their Theorem~4; their second case
should be solving $\lambda(ps-qr)-npr=0$, without an $n$ in the exponent.)

In the first case, the proof here follows the proof in \bb.

The second case is even easier. If \hbox{$p=0$,} then the condition simplifies to \hbox{$qr=0$,} and if \hbox{$r=0$,}
it simplifies to \hbox{$ps=0$;} in both cases, $n$ is arbitrary.
Otherwise, it is a matter of computing $n={\lambda(p s-qr)\over p r}$ and determining whether it
is actually a nonnegative integer. This ends Proof~\#2.

\medskip

The current author discovered a completely new proof, not considering the Jordan form of $A$. 
To handle one case, a lemma is required:

\newlemma\cosine{{\bf(Lemma~6 \bb, originally from \shank.)} 
\item{\rm1.} The following decision problem is decidable.
\itemitem{} Instance: rational numbers $p,q\in[-1,1]$.
\itemitem{} Question: Does there exist $\theta\in\R$ and an integer $n\in\N$ with $\cos(\theta)=p$ and
	$\cos(n\theta)=q$?
\item{\rm2.} When $p\not\in\set{0, {1\over2},1}$ there is a finite number of such $n$, and those values
	can be computed effectively.
}

\medskip

If $A$ and $B$ are matrices, the notation $A\sim B$ will be used when $A$ is a nonzero multiple of $B$. This is an
equivalence relation that preserves matrix multiplication. (If $A\sim B$ and $C\sim D$, then $AC\sim BD$.)

\medskip

\vbox{
\newlemma\power{The following decision problem is decidable. A value of $n$ is effectively computable as well,
	in the case of a YES answer.
\item{} Instance: A rational $2\times2$ matrix $A$.
\item{} Question: Does there exist a positive integer $n$ such that $A^n \sim I$?
}
\par}

\medskip

\ital{Proof of \power.} Suppose that $A^k=r\,I$, and 
write $A=PJP^{-1}$ where $J$ is in Jordan form. Then $A^k= PJ^k P^{-1}$, and
$$J^k = P^{-1} A^k P = P^{-1} r\,I P = r\,I.$$
This implies that $J$ is itself a diagonal matrix, with the eigenvalues of $A$ on its diagonal,
and $A$ is diagonalizable.
Thus, for any eigenvalue $\lambda$ of $A$, $\lambda^k=r$, which implies that $\lambda^k = \mu^k$ if $\lambda,\mu$ are
eigenvalues of $A$, and ${\lambda\over\mu}$ is a root of unity.

Now, we provide an algorithm for establishing whether $A^k\sim I$ for a general matrix $A$. We
assume that $A$ is invertible and that $A \not\sim I$, which are easily seen to be decidable. 
Otherwise, compute the eigenvalues of $A$ (the roots of $\lambda^2+b\lambda+c$) and check for diagonalizability.
If $A$ is not diagonalizable, then $A^k \not\sim I$ for any $k$. 
If the eigenvalues of $A$ are equal, then $A\sim I$ (which we have eliminated already).

If $\lambda_1$ and $\lambda_2$ are the eigenvalues of $A$, then write $P=\brcmat{\alpha&\beta\cr \gamma&\delta\cr}$ and 
$J=\brcmat{\lambda_1&0\cr 0&\lambda_2\cr}$, so that $A=PJP^{-1}$. A bit of computation shows that
$$A^k = \brcmat{\alpha\delta \cdot e^k - \beta\gamma \cdot f^k & -\alpha \beta \cdot e^k + \alpha \beta \cdot f^k \cr
	\gamma\delta \cdot e^k - \gamma\delta \cdot f^k & -\beta\gamma \cdot e^k + \alpha \delta \cdot f^k\cr}.$$

In order to have $A^k\sim I$, we must have $\alpha \beta(\lambda_1^k-\lambda_2^k)=0$ (to make the
$(1,2)$ entry be zero), $\gamma \delta(\lambda_1^k - \lambda_2^k)=0$ (for the $(2,1)$ entry), and
$(\lambda_1^k - \lambda_2^k)(\alpha \delta + \beta \gamma)=0$ (after equating the $(1,1)$ and $(2,2)$ entries).
If $\lambda_1^k \not=\lambda_2^k$, then the solutions to these equations have $\alpha = \gamma=0$ or $\beta = \delta =0$;
in either case, $P$ is not invertible. Hence, $\lambda_1^k = \lambda_2^k$, which means $\rho={\lambda_1 \over \lambda_2}$
is a root of unity.

Computation shows that $\rho ={b^2-2c\over 2c} + {b\over 2c}\sqrt{b^2-4c}$. 
If $b^2-4c \ge 0$, then we check to see whether $\rho=\pm1$, as these are the only real roots of unity. 
This can be done by checking whether
$b=0$ (making $\rho=-1$) or $b^2=4c$ (making the $\rho=1$). 
If both conditions are false, then $A^k \not\sim I$ for any $k$; if either is true, then we can take $k=2$.

If $b^2-4c < 0$, then $\rho={b^2-2c\over 2c} + {b\over 2c}\sqrt{4c-b^2}\cdot i$, which is
a complex number with both coefficients being real. In order to have $\rho$ be a root of unity, we must first
have $\rho \overline\rho=1$, which is equivalent to the condition
$$-b^2(c+1)(c-1)(c^2+1)(b^2-4c)=0,$$
where we only need check the first three factors. If all three are nonzero, then $A^k \not\sim I$.

Now, since $\size\rho =1$, we can write $\rho = \cos\theta + i\sin \theta$ for some angle $\theta$. Since
$\rho$ is a root of unity, we seek a positive integral solution to $1 = \rho^n = \cos(n\theta)+i\sin(n\theta)$.
This can be done by applying \cosine\ with $p= {b^2-2c \over 2c}$ and $q=1$.

\medskip

\ital{Proof \#3 of \one.} We need to determine whether there are any nonnegative integer
solutions to the equation $B A^k B = 0$. Assume that $B\not=0$.
Then let the characteristic polynomial of $A$ be \hbox{$\lambda^2 + b \lambda +c=0$}.
(Note that $b,c$ are rational, and integral if $A$ is.)
The Cayley-Hamilton Theorem \chthm\ implies that \hbox{$A^2 + b A + c I=0$}.
Define $r_k$ by: 
$$r_k = \cases{0 & if $k=1$\cr {c \over b -r _{k-1}}& otherwise.\cr}$$
An easy induction argument shows that $A^k \sim A + r_k I$, for all $k\ge1$. 

It is possible for $r_k$ to be undefined for some value of $k$. This will happen only if $r_{k-1}=b$, i.e., if
$A^{k-1} \sim A + bI$. Thereupon,
$$A^k = A \cdot A^{k-1} = A(A + bI) = A^2+bA = -c I \sim I.$$
We eliminate this possibility by use of \power: We check to see whether $A^k \sim I$ for some positive integer $k$.
If so, then $\set{A,B}$ is mortal iff $0\in \set{BA^iB \mid 0 \le i < k}$.

Now, substitution into the equation
$BA^k B=0$ and manipulation results in the equation $-BAB = r_k\cdot B^2$.

There is at most one real number that makes this equation true. If $-BAB$ is not a multiple of $B^2$,
then $\set{A,B}$ is not mortal, and the algorithm finishes. Otherwise, let $x$ be this multiple. Then the equation
$r_k=x$ can be solved for $k$; in fact, there is a formula that determines the unique value of $k$:
$$
k = {\dst\ln\({-\sqrt{b ^{2} -4c }-b  + 2x \over \sqrt{b ^{2} -4c }-b  + 2x}\) 
\over 
\dst-\ln\({-c  \over b  + \sqrt{b ^{2} -4c }}\) + \ln\({c  \over -b  + \sqrt{b ^{2} -4c }}\)}
$$
(Exact evaluation is unnecessary;
if this value is determined to be within 0.1 of an integer $k_0$, we then check whether $BA^{k_0} B=0$.)

This concludes Proof~\#3.

\newsection\mmInv{Matrix Mortality for Other Numbers of Invertible Matrices}

This method can be used to establish relationships between the decidability of the Matrix Mortality Problem depending
on the number invertible and the number of noninvertible matrices.

Given a fine set $S$ of non-invertible (singular) matrices and a finite set $I$ of invertible matrices, let 
$M(S,I)$ denote the proposition that $S \cup I$ is mortal. Also let $D(s,i)$ denote the proposition
that the MMP is decidable, given any set $S$ of non-invertible matrices with $\size S = s$ and any set $I$ of invertible
matrices with $\size I = i$.

Note that if $T$ is any set of matrices, its mortality is equivalent to the mortality of a set obtained from $T$ by
replacing (or by adding) an element $M$ of $T$ with an element $M'$ such that $M\sim M'$; this is a consequence of
\product. This fact will be used implicitly in what follows.

\newthm\two{For all positive integers $i,j,k$, $D(i,k)$ iff $D(j,k)$; that is, decidability does not depend on the
number of non-invertible matrices.}

\vbox{
\ital{Proof of \two.} The theorem will be proven for the case $i=2$, and transitivity guarantees that the rest of the cases
will follow.

Assume that $D(2,k)$ is true, that there is an algorithm that establishes whether $M(S,I)$ is true when
$\size S=2$ and $\size I=k$.

First, suppose $j=1$; that is, $S$ contains one non-invertible matrix $B$. An algorithm to determine 
whether $M(\set B, I)$ is as follows: If $B=0$, then the answer is yes. Otherwise, determine whether $M(\set{B,-B},I)$
is true. Thus, $D(2,k)\to D(1,k)$.

}

Now suppose that $j>2$. To establish whether $M(S,I)$ is true, we check for the equivalent
condition given by \product, where at most two non-invertible matrices appear. Thus, we simply determine whether
$M(\set{B,-B},I)$ is true for every $B\in S$, and whether $M(\set{B_1,B_2},I)$ is true for every $B_1,B_2\in S$
with $B_1\not=B_2$. Our assumed algorithm will be repeatedly run on these instances, and return true if any one
of the cases is true; otherwise, it returns false.

Now, suppose $D(i,k)$ is true, that there is an algorithm that determines whether $M(S,I)$ is true, where
$\size S=i$ and $\size S=k$. 

If $i>2$, then we can find a set $S' \supseteq S$ which contains only nonzero multiples of elements of $S$,
and such that $\size {S'}=i$; then
we run our algorithm on $S'$ and $I$ to determine whether $M(S',I)$ is true (and hence also $M(S,I)$).

The last case, where $i=1$, is the trickiest; it is actually a re-working of Proof \#1 of \bb above. 
Let $S=\set{B_1,B_2}$. If $0\in S$, then $M(S,I)$ is trivially true.
Otherwise, write $B_1=a b^\top$ and $B_2=c d^\top$ (as can be done for rank~one $2\times2$ matrices).
Now, \bb\ implies that $M(2,k)$ will be true if
there is a product $B_{n_1} A_1 \cdots A_n B_{n_2}=0$ such that $A_i\in I$, for $1\le i\le n$, and $n_1,n_2\in\set{1,2}$.
To test whether this is the case, we simply check every single case: $(n_1,n_2)=(1,1)$ or $(2,2)$ can be checked
directly by the assumed algorithm; $(n_1,n_2)=(1,2)$ can be checked by determining whether $M(\set{B'}, I)$ is true,
where $B'=c b^\top$.

Note that $M(\set{B'},I)$ is true iff:
\smallskip
\bii there is a product $B' A_1 \cdots A_n B'=0$ such that $n\ge0$ and $A_i \in I$; which is true, iff
\bii there is a product $c b^\top A_1 \cdots A_n c b^\top=0$ such that [\dots]; which is true, iff
\bii there is a product $b^\top A_1 \cdots A_n c=0$ such that [\dots]; which is true, iff
\bii there is a product $a b^\top A_1 \cdots A_n c d^\top=0$ such that [\dots]; which is true, iff
\bii there is a product $B_1 A_1 \cdots A_n B_2=0$ such that [\dots],
\smallskip
\noindent which is equivalent to one of the cases that we want to check, namely, $(n_1,n_2)=(1,2)$.

Similarly, $M(\set{a d^\top}, I)$ will establish whether another case --- $(n_1,n_2)=(2,1)$ --- is true.

This establishes the truth of \two.

\newsection\refs{References}

\newcount\refnum \refnum=0
\def\nextref#1{\global\advance\refnum by 1 \smallskip
	\immediate\write3{\def\string#1{\hbox{\noexpand\rm[\the\refnum]}}}
	\item{#1}}

\immediate\openout3=mortalityrefs.aux

\nextref\bb O.~Bournez and M.~Branicky. The Morality Problem for Matrices of Low Dimensions, \ital{Theory of Computing
Systems} 35, 433--448 (2002).

\nextref\twentyeighteen Julien Cassaigne, Vesa Halava, Tero Harju, Fran\c cois Nicolas.
	Tighter Undecidability Bounds for Matrix Mortality, Zero-in-the-Corner Problems, and More.
	{\tt arXiv:1404.0644v3 [cs.DM] 5 Sep 2014 }

\nextref\chthm A.~Cayley. A Memoir on the Theory of Matrices. \ital{Philosophical Transactions} 148 (1858). 17--37.
	(Historical note: Cayley proved the $2\times2$ and $3\times3$ cases in this paper, which is what we need here. The
	general case was proven in 1878 by G.~Frobenius.)

\nextref\nucrod C.~Nuccio and E.~Rodardo. Mortality problem for $2 \times 2$ integer matrices. In \ital{Lecture
Notes in Computer Science}, volume 4910 of \ital{Proceeding SOFSEM'08 Proceedings of
the 34th conference on Current trends in theory and practice of computer science},
pp.~400--405. Springer-Verlag, 2008.

\nextref\pater M.~S.~Paterson. Unsolvability in $3\times3$ matrices. \ital{Studies in Applied Mathematics},
XLIX(1):105--107, March 1970.

\nextref\shank H.~S.~Shank. The rational case of a matrix problem of Harrison. \ital{Discrete Mathematics}, 28:207--212,
1979.

\nextref\twentyfive N.~Vereshchagin. Occurrence of zero in a linear recursive sequence. \ital{Mathematical Notes}, 
38:609--615, 1985. Translation from \ital{Matematicheskie Zametki} 38(2):609--615, 1985.

\closeout3

% Developments in language theory.QA267.3 .C63 2009  

\vfill\vfill\vfill\vfill\vfill\vfill
\line{\hfill This document last modified \lastmod.}
\closeout1
\eject

\end %%%%%% 

\advance\vsize by 0.25in

\ignore{
\noindent *** Not part of the proof, but looking to establish \one\ using yet another approach ***
\medskip

If $k\ge 1$, then we can write $A^k=a_k A + b_k I = a_k\(A + {b_k\over a_k}\, I\)\sim A + {b_k\over a_k}\, I$, 
for some integers $a_k,b_k$, due to the Cauchy-Hamilton Theorem and polynomial division. Suppose that
$A^2+\alpha A+\beta I =0$ (i.e., the characteristic polynomial of $A$ is $\lambda^2 +\alpha\lambda + \beta$).

The base case is $(a_1,b_1)=(1,0)$. If $k\ge1$, then
\medskip
\centerline{$\eqalign{
a_{k+1} A + b_{k+1} I &= A^{k+1} = A \cdot A^k 
= A(a_k A + b_k I) = a_k A^2 + b_k A
a_k (-\alpha A-\beta I) + b_k A
= (b_k - a_k \alpha)A + (-a_k \beta)I.\cr
}$}
\medskip
This is an order one linear recurrence for $(a_k,b_k)$.

$$r_{k+1}\equiv {b_{k+1} \over a_{k+1}} = {-a_k \beta \over b_k - a_k \alpha}
= {-a_k \beta /a_k \over (b_k - a_k \alpha)/a_k}
= {-\beta \over\dst {b_k \over a_k} -\alpha}
= {-\beta \over r_k - \alpha}
$$
Iterating this produces formulas of the form
$$\alpha - {\beta \over \dst \alpha - {\beta \over \dst \alpha - {\beta \over \dst \alpha - {\beta \over \dst \alpha}}}}$$
which, when simplified, appear as follows. 
$$
\brcmat{1 & 0\cr\cr 2 & {\beta  \over \alpha }\cr\cr 3 & {\beta \alpha  \over \alpha ^{2}-\beta }\cr\cr 4 & {\beta \(-1\alpha ^{2} + \beta \) \over -1\alpha ^{3} + 2\alpha \beta }\cr\cr 5 & {\beta \alpha \(\alpha ^{2} -2\beta \) \over \alpha ^{4} -3\alpha ^{2}\beta  + \beta ^{2}}\cr\cr 6 & {\beta \(\alpha ^{4} -3\alpha ^{2}\beta  + \beta ^{2}\) \over \alpha ^{5} -4\alpha ^{3}\beta  + 3\alpha \beta ^{2}}\cr\cr 7 & {-3\beta \alpha \(-1\alpha ^{2} + \beta \)\({-1\alpha ^{2} \over 3} + \beta \) \over -1\alpha ^{6} + 5\alpha ^{4}\beta  -6\alpha ^{2}\beta ^{2} + \beta ^{3}}\cr\cr 8 & {\beta \(-1\alpha ^{6} + 5\alpha ^{4}\beta  -6\alpha ^{2}\beta ^{2} + \beta ^{3}\) \over -1\alpha ^{7} + 6\alpha ^{5}\beta  -10\alpha ^{3}\beta ^{2} + 4\alpha \beta ^{3}}\cr\cr 9 & {\beta \alpha \(\alpha ^{6} -6\alpha ^{4}\beta  + 10\alpha ^{2}\beta ^{2} -4\beta ^{3}\) \over \alpha ^{8} -7\alpha ^{6}\beta  + 15\alpha ^{4}\beta ^{2} -10\alpha ^{2}\beta ^{3} + \beta ^{4}}\cr\cr 10 & {\beta \(\alpha ^{8} -7\alpha ^{6}\beta  + 15\alpha ^{4}\beta ^{2} -10\alpha ^{2}\beta ^{3} + \beta ^{4}\) \over \alpha \(\alpha ^{8} -8\alpha ^{6}\beta  + 21\alpha ^{4}\beta ^{2} -20\alpha ^{2}\beta ^{3} + 5\beta ^{4}\)}\cr\cr 11 & {\beta \alpha \(\alpha ^{8} -8\alpha ^{6}\beta  + 21\alpha ^{4}\beta ^{2} -20\alpha ^{2}\beta ^{3} + 5\beta ^{4}\) \over \alpha ^{10} -9\alpha ^{8}\beta  + 28\alpha ^{6}\beta ^{2} -35\alpha ^{4}\beta ^{3} + 15\alpha ^{2}\beta ^{4}-\beta ^{5}}\cr\cr 12 & {\beta \(-1\alpha ^{10} + 9\alpha ^{8}\beta  -28\alpha ^{6}\beta ^{2} + 35\alpha ^{4}\beta ^{3} -15\alpha ^{2}\beta ^{4} + \beta ^{5}\) \over -1\alpha ^{11} + 10\alpha ^{9}\beta  -36\alpha ^{7}\beta ^{2} + 56\alpha ^{5}\beta ^{3} -35\alpha ^{3}\beta ^{4} + 6\alpha \beta ^{5}}\cr}
$$
(Note that factoring doesn't appear to help;
$$r_{11} = {\beta \alpha \(\alpha ^{4} -3\alpha ^{2}\beta  + \beta ^{2}\)\(\alpha ^{4} -5\alpha ^{2}\beta  + 5\beta ^{2}\) \over \alpha ^{10} -9\alpha ^{8}\beta  + 28\alpha ^{6}\beta ^{2} -35\alpha ^{4}\beta ^{3} + 15\alpha ^{2}\beta ^{4}-\beta ^{5}}.)$$

Odd terms:
$$
\brcmat{1 & 0\cr\cr 3 & {\beta \alpha  \over \alpha ^{2}-\beta }\cr\cr 5 & {\beta \alpha ^{3} -2\beta ^{2}\alpha  \over \alpha ^{4} -3\alpha ^{2}\beta  + \beta ^{2}}\cr\cr 7 & {\beta \alpha ^{5} -4\beta ^{2}\alpha ^{3} + 3\beta ^{3}\alpha  \over \alpha ^{6} -5\alpha ^{4}\beta  + 6\beta ^{2}\alpha ^{2}-\beta ^{3}}\cr\cr 9 & {\beta \alpha ^{7} -6\beta ^{2}\alpha ^{5} + 10\beta ^{3}\alpha ^{3} -4\beta ^{4}\alpha  \over \alpha ^{8} -7\alpha ^{6}\beta  + 15\alpha ^{4}\beta ^{2} -10\alpha ^{2}\beta ^{3} + \beta ^{4}}\cr\cr 11 & {\beta \alpha ^{9} -8\beta ^{2}\alpha ^{7} + 21\beta ^{3}\alpha ^{5} -20\beta ^{4}\alpha ^{3} + 5\beta ^{5}\alpha  \over \alpha ^{10} -9\alpha ^{8}\beta  + 28\alpha ^{6}\beta ^{2} -35\alpha ^{4}\beta ^{3} + 15\alpha ^{2}\beta ^{4}-\beta ^{5}}\cr\cr 13 & {\beta \alpha ^{11} -10\beta ^{2}\alpha ^{9} + 36\beta ^{3}\alpha ^{7} -56\beta ^{4}\alpha ^{5} + 35\beta ^{5}\alpha ^{3} -6\beta ^{6}\alpha  \over \alpha ^{12} -11\alpha ^{10}\beta  + 45\alpha ^{8}\beta ^{2} -84\alpha ^{6}\beta ^{3} + 70\alpha ^{4}\beta ^{4} -21\alpha ^{2}\beta ^{5} + \beta ^{6}}\cr\cr 15 & {\beta \alpha ^{13} -12\beta ^{2}\alpha ^{11} + 55\beta ^{3}\alpha ^{9} -120\beta ^{4}\alpha ^{7} + 126\beta ^{5}\alpha ^{5} -56\beta ^{6}\alpha ^{3} + 7\beta ^{7}\alpha  \over \alpha ^{14} -13\alpha ^{12}\beta  + 66\alpha ^{10}\beta ^{2} -165\alpha ^{8}\beta ^{3} + 210\alpha ^{6}\beta ^{4} -126\alpha ^{4}\beta ^{5} + 28\alpha ^{2}\beta ^{6}-\beta ^{7}}\cr}
$$

Even terms:
$$
\brcmat{2 & {\beta  \over \alpha }\cr\cr 4 & {\alpha ^{2}\beta -\beta ^{2} \over \alpha ^{3} -2\beta \alpha }\cr\cr 6 & {\alpha ^{4}\beta  -3\beta ^{2}\alpha ^{2} + \beta ^{3} \over \alpha ^{5} -4\beta \alpha ^{3} + 3\beta ^{2}\alpha }\cr\cr 8 & {\alpha ^{6}\beta  -5\alpha ^{4}\beta ^{2} + 6\alpha ^{2}\beta ^{3}-\beta ^{4} \over \alpha ^{7} -6\beta \alpha ^{5} + 10\beta ^{2}\alpha ^{3} -4\beta ^{3}\alpha }\cr\cr 10 & {\alpha ^{8}\beta  -7\alpha ^{6}\beta ^{2} + 15\alpha ^{4}\beta ^{3} -10\alpha ^{2}\beta ^{4} + \beta ^{5} \over \alpha ^{9} -8\beta \alpha ^{7} + 21\beta ^{2}\alpha ^{5} -20\beta ^{3}\alpha ^{3} + 5\beta ^{4}\alpha }\cr\cr 12 & {\alpha ^{10}\beta  -9\alpha ^{8}\beta ^{2} + 28\alpha ^{6}\beta ^{3} -35\alpha ^{4}\beta ^{4} + 15\alpha ^{2}\beta ^{5}-\beta ^{6} \over \alpha ^{11} -10\beta \alpha ^{9} + 36\beta ^{2}\alpha ^{7} -56\beta ^{3}\alpha ^{5} + 35\beta ^{4}\alpha ^{3} -6\beta ^{5}\alpha }\cr\cr 14 & {\alpha ^{12}\beta  -11\alpha ^{10}\beta ^{2} + 45\alpha ^{8}\beta ^{3} -84\alpha ^{6}\beta ^{4} + 70\alpha ^{4}\beta ^{5} -21\alpha ^{2}\beta ^{6} + \beta ^{7} \over \alpha ^{13} -12\beta \alpha ^{11} + 55\beta ^{2}\alpha ^{9} -120\beta ^{3}\alpha ^{7} + 126\beta ^{4}\alpha ^{5} -56\beta ^{5}\alpha ^{3} + 7\beta ^{6}\alpha }\cr\cr 16 & {\beta \alpha ^{14} -13\beta ^{2}\alpha ^{12} + 66\beta ^{3}\alpha ^{10} -165\beta ^{4}\alpha ^{8} + 210\beta ^{5}\alpha ^{6} -126\beta ^{6}\alpha ^{4} + 28\beta ^{7}\alpha ^{2}-\beta ^{8} \over \alpha ^{15} -14\beta \alpha ^{13} + 78\beta ^{2}\alpha ^{11} -220\beta ^{3}\alpha ^{9} + 330\beta ^{4}\alpha ^{7} -252\beta ^{5}\alpha ^{5} + 84\beta ^{6}\alpha ^{3} -8\beta ^{7}\alpha }\cr}
$$

$$\lim_{k\to\infty} r_k = {\alpha  \pm \sqrt{\alpha ^{2} -4\beta } \over 2}
$$

Now, let's let Maple solve $r_{k+1}(\alpha-r_k) = \beta$ with $r(1)=0$:
$$
r_k = 
{\(-\sqrt{\alpha ^{2} -4\beta } + \alpha \)\(2\beta  \over\dst -\alpha  + \sqrt{\alpha ^{2} -4\beta }\)^{k} -\(\alpha  + \sqrt{\alpha ^{2} -4\beta }\)\(-2\beta  \over\dst \alpha  + \sqrt{\alpha ^{2} -4\beta }\)^{k} \over 
\dst 2 \cdot \(2\beta  \over -\alpha  + \sqrt{\alpha ^{2} -4\beta }\)^{k} -2\(-2\beta  \over \alpha  + \sqrt{\alpha ^{2} -4\beta }\)^{k}}
$$
Amazing. Double amazing is the fact that we can solve $r_k=x$ for $k$:
$$
k = {\dst\ln\({-\sqrt{\alpha ^{2} -4\beta }-\alpha  + 2x \over \sqrt{\alpha ^{2} -4\beta }-\alpha  + 2x}\) 
\over
\dst \ln\({\alpha  - \sqrt{\alpha ^{2} -4\beta }  \over \alpha  + \sqrt{\alpha ^{2} -4\beta }}\)}
$$
unsimplified form:
$$
{\dst\ln\({-\sqrt{\alpha ^{2} -4\beta }-\alpha  + 2x \over \sqrt{\alpha ^{2} -4\beta }-\alpha  + 2x}\) 
\over 
\dst-\ln\({-\beta  \over \alpha  + \sqrt{\alpha ^{2} -4\beta }}\) + \ln\({\beta  \over -\alpha  + \sqrt{\alpha ^{2} -4\beta }}\)}
$$

This means there's a case-free proof that with only one invertible matrix, the problem is decidable.

}

\ignore{
\noindent *** Now to reconstruct what I did in the doctor's office, Aug 30. ***

Making the assumption that $B_1, B_2$ are not invertible, $A_1,A_2$ are. Write
$A_i = P_i J_i P_i^{-1}$, 
$Q = P_1^{-1} P_2 = \brcmat{p&q\cr r&s\cr}$. 
Then $A_i^{\,k} = P_i J_i^{\,k} P_i^{-1}$ and
now we attempt to find a nice form for $B_1 \prod_{i=1}^n A_{n_i} B_2$. Group the product by repeat values of
$A_i$, to get $B_1 A_1^{\,k} A_2^{\,\ell} \cdots A_?^{\,?} B_2$.
One run of $A_1$'s has been done (B\&B); now we try two and more. 

Write $B_1 P_1 = * b^\top$, $P_2^{-1} B_2 =  a *^\top$, $b^\top = [t~u]$, $a = \brcmat{v\cr w\cr}$; 
then the matrix equation with 2 runs is
$$\eqalignno{
0 = B_1 A_1^{\,k}A_2^{\,\ell}B_2 &= B_1 P_1 J_1^{\,k} P_1^{-1} P_2 J_2^{\,\ell} P_2^{-1} B_2
= (B_1 P_1) J_1^{\,k} (P_1^{-1} P_2) J_2^{\,\ell} (P_2^{-1} B_2)= *b^\top J_1^{\,k} Q J_2^{\,\ell} a*^\top\cr
0&=[t~u] \cdot J_1^{\,k} \brcmat{p&q\cr r&s\cr} J_2^{\,\ell} \brcmat{v\cr w\cr}&(2)
\cr
}$$

Write $B_1 P_1 = * b^\top$, $P_1^{-1} B_2 =  a *^\top$ --- Yes, this is confusing, because it's not the same $a$ as
in (2), but ``there aren't enough good letters'' (RT) --- $b^\top = [t~u]$, $a = \brcmat{v\cr w\cr}$; 
the matrix equation with 3 runs is
$$\eqalignno{
0 = B_1 A_1^{\,k}A_2^{\,\ell}A_1^{\,m} B_2 &= B_1 P_1 J_1^{\,k} P_1^{-1} P_2 J_2^{\,\ell} P_2^{-1}P_1 J_1^{\,m} B_2
= (B_1 P_1) J_1^{\,k} (P_1^{-1} P_2) J_2^{\,\ell} (P_2 P_1^{-1}) J_1^{\,m} (P_1^{-1} B_2)
\cr&\sim *b^\top J_1^{\,k} Q J_2^{\,\ell} (\adj Q) J_1^{\,m} a*^\top\cr
0&=[t~u] \cdot J_1^{\,k} \brcmat{p&q\cr r&s\cr} J_2^{\,\ell} \brcmat{s&-q\cr -r&p\cr}  J_1^{\,m} \brcmat{v\cr w\cr}
&(3)
\cr
}$$
\medskip

\noindent
Now, it breaks into 4 cases. 
\smallskip
\itemitem{(a)} The easy one is: neither $A_i$ is diagonalizable.
\itemitem{(b)} Exactly one is, the first.
\itemitem{(c)} Exactly one is, the second.
\itemitem{(d)} Both are diagonalizable.
\medskip

\noindent\llap{Case (a): }
Write
$J_1^{\,k}\sim \brcmat{L&k\cr 0&L\cr}$ ($\sim$ means ``is a nonzero multiple of'') and 
$J_2^{\,\ell} \sim \brcmat{M&\ell\cr 0&M\cr}$.

\noindent \llap{Subcase (a2): }
Equation (2) becomes
$$
[rtw]k\ell + [Lptw + Lruw] \ell + [Mrtv + Mstw]k + [LMptv + LMqtw + LMruv + LMsuw]=0
$$
which ought to be solvable/decidable. (Quantities inside of brackets are known to us.)\footnote*{
\dots\ and is, according to {\tt https://mathoverflow.net/questions/\-207482/\-algorithmic-un-solvability-of-diophantine-equations-of-given-degree-with-given}, the solvability for the polynomial equation
\hbox{$P(x_1,\ldots,x_k)=0$}  could depend on $d$,
the total degree of $P$ ($k$ is the number of variables here). 
\smallskip
\bi Clearly, $(d,1),(1,k)$ are decidable. 
\bi $(2,k)$ is decidable (Siegel, 1972).
\bi $(d,k)=(3,2)$ is decidable (Baker and Coates, ``Integer Points on Curves of Degree 1'', 1970, and Siegel).
\medskip
\bi $(3,3)$ is undecidable; this would also allow the ranks of elliptic curves to be determined (which is also
	unsolved).
\bi $(d,k)$ with $d\ge 4$ is NOT decidable if $k$ is ``sufficiently large'' (Thoralf Skolem, 1920s), because
	any Diophantine equation is equivalent to one with degree $\le4$.
\bi $(d,k)$ with $k\ge 11$ is NOT decidable if $d$ is ``sufficiently large''.
	(Zhi Wei Sun, 1992, ``Reduction of unknowns in Diophantine representations'')
\bi If $k\ge 9$, then the equation is NOT decidable. 9 is not known to be the largest upper bound.
\medskip
\noindent If the solutions are supposed to be positive:
\bi $(4,58),(8,38),(12,32),(16,29),(20,28),(24,26),(28,25),(36,24),(96,21),(2668,19),(2\cdot 10^5, 14)$,
	\hbox{$(6.6\cdot 10^{43},13)$,} $(1.3\cdot 10^{44}, 12), (4.6\cdot 10^{44}, 11), (8.6\cdot10^{44},10), 
	(1.6\cdot 10^{45},9)$ are 
	undecidable (Jones, ``Universal Diophantine equation'', 1982).
\par}

\medskip
\noindent\llap{Subcase (a3): } Equation (3) becomes
$$
\matrix{\hfilneg0 = [-r^{2}tw]k \ell m
	+ 
[L prtw -Lr^{2}tv ]k \ell
	+ 
[ -L prtw -L r^{2}uw]\ell m
	+ 
[LMpstw -LMqrtw ]k 
	\cr
\quad\quad+
[L^{2} p^{2}tw -L^{2} prtv + L^{2} pruw -L^{2} r^{2}uv  ]\ell
	+
[LMpstw -LMqrtw ]m\hfilneg\cr
\quad\quad+ 
[L^{2}Mpstv + L^{2}Mpsuw -L^{2}Mqrtv -L^{2}Mqruw   ]\hfilneg\cr}
$$
which is degree 3 ($k\ell m$).

\medskip
\noindent\llap{Subcase (a4): } Equation (4) becomes a polynomial with total degree 4 (there is a $k\ell mn$ term):
\medskip
\centerline{$
\eqalign{
0=k\ell mn&~[-r^3tw]\cr
{}+k\ell m&~[-Mr^2stw-Mr^3tv]\cr
{}+\ell mn&~[-Lpr^2tw - Lr^3uw]\cr
{}+k\ell&~[LMprstw-LMqr^2tw]\cr
{}+kn&~[LMprstw-LMqr^2tw]\cr
{}+\ell m&~[-LMprstw-LMr^2suw-LMpr^2tv-LMr^3uv]\cr
{}+mn&~[LMprstw-LMqr^2tw]\cr
{}+k&~[LM^2prstv-LM^2qrstw+LM^2ps^2tw-LM^2qr^2tv]\cr
{}+\ell&~[-L^2Mpqrtw+L^2Mprsuw+L^2Mp^2stw-L^2Mqr^2uw]\cr
{}+m&~[LM^2prstv-LM^2qrstw+LM^2ps^2tw-LM^2qr^2tv]\cr
{}+n&~[-L^2Mpqrtw+L^2Mprsuw+L^2Mp^2stw-L^2Mqr^2uw]\cr
{}+&~[L^{2}M^{2}p^{2}stv -L^{2}M^{2}pqrtv + L^{2}M^{2}pqstw + L^{2}M^{2}prsuv + L^{2}M^{2}ps^{2}uw -L^{2}M^{2}q^{2}rtw -L^{2}M^{2}qr^{2}uv -L^{2}M^{2}qrsuw]\cr
}
$}

\medskip
\noindent \llap{Subcase (b2): } OTOH, write $J_1^{\,k}=\brcmat{\lambda^k &0\cr 0&\mu^k\cr}$ and
$J_2^{\,\ell} \sim \brcmat{L&\ell\cr 0&L\cr}$. Now, (2) becomes
$$
\lambda ^{k}[Lptv + Lqtw] + \mu ^{k}[Lruv + Lsuw] + \lambda ^{k}\ell[ptw] + \mu ^{k}\ell[ruw] = 0
$$
(Same form, just variables shuffled, and $k$ and $\ell$ swapped.)

\medskip
\noindent \llap{Subcase (c2): }Write $J_1^{\,k} \sim \brcmat{L&k\cr 0&L\cr}$ and 
$J_2^{\,\ell}=\brcmat{\lambda^\ell &0\cr 0&\mu^\ell\cr}$. Now, (2) becomes
$$
\lambda ^{\ell}[Lptv+Lruv] + \mu ^{\ell}[Lqtw+Lsuw] + \lambda ^{\ell}k[rtv] + \mu ^{\ell}k[stw]=0
$$

\medskip
\noindent \llap{Subcase (d2): }Write $J_1^{\,k} = \brcmat{L^k&0\cr 0&M^k\cr}$ and 
$J_2^{\,\ell}=\brcmat{\lambda^\ell &0\cr 0&\mu^\ell\cr}$. Now, (2) becomes
$$
[ptv]\lambda ^{\ell}L^{k} + [ruv]\lambda ^{\ell}M^{k} + [qtw]\mu ^{\ell}L^{k} + [suw]\mu ^{\ell}M^{k}
= 0
$$

\medskip
\noindent \llap{Subcase (d3): }\medskip
\centerline{$
[pstv]L^{m + k}\lambda ^{\ell} 
+[-qrtv]L^{m+k}\mu^\ell 
+ [-qruw]M^{m + k} \lambda ^{\ell}
+ [psuw]M^{m+k}\mu ^{\ell} 
+ [pqtw] \(\mu ^{\ell}-\lambda ^{\ell}\) L^k M^m
+ [rsuv]\(\mu ^{\ell}-\lambda ^{\ell}\) L^{m}M^{k}
= 0
$}

\medskip
\noindent\llap{Subcase (d4) } \dots\ is an ungodly mess.
}

*** Double Schur Triangulation.

Given two matrices $A_1,A_2$, there exist unitary (orthogonal) matrices $S$ and $T$ such that
$A_1 = ST_1 T$, $A_2 = ST_2 T$, and $T_1,T_2$ are both upper-triangular.

Furthermore, we may assume that $T_i=\brcmat{1&a_i\cr 0&b_i\cr}$ for $i=1,2$, if $A_1,A_2$ are invertible
$2\times2$ matrices. Then
$B_1A_*A_*\cdots A_*B_2 = 0$ iff
$$\eqalign{
B_1\cdot S T_* T \cdot S T_* T \cdots S T_* T \cdot B_2  &=\bf0\cr
(B_1S) T_* (TS) T_* (TS) \cdots T_* (TB_2) &=\bf0\cr
}
$$
Write $B_1 S=*\brcmat{g~h}$, $B_2=S\brcmat{c&d\cr e&f\cr}^{-1}\brcmat{i\cr j\cr}*$, and $TS = \brcmat{c&d\cr e&f\cr}$. Then
$$\brcmat{g~h} \cdot T_*\brcmat{c&d\cr e&f\cr} \cdot T_* \brcmat{c&d\cr e&f\cr} \cdots
	T_* \brcmat{c&d\cr e&f\cr} \cdot \brcmat{c&d\cr e&f\cr}^{-1}\cdot\brcmat{i\cr j\cr}=0.$$

The LHS has the form
$$L_k = \brcmat{g~h} \cdot \(T_* \brcmat{c&d\cr e&f\cr}\)^k \brcmat{c&d\cr e&f\cr}^{-1}\brcmat{i\cr j\cr},$$
and there might be a recurrence that connects $L_k$ to $L_{k+1}$. (The inverse of $TS$ is here to simplify the
calculations.)

This raises the following question about decidability of the following problem:
Given two recurrences $a_k = f_1(a_{k-1})$, $a_k = f_2(a_{k-1})$ and $a_0$, is it possible to choose which recurrence
is applied at each step, to end up with $a_n=0$ for some $n$?

Getting back to the problem of formulating the recurrence, note that $a_n$ really means $a_{i_n}$, where
$i_n\in\set{1,2}$, same for $b_n$. 

$$
\eqalign{
L_0 &= hc -gd\cr\cr
L_1 &=  gja_{1} + hjb_{1} + g\cr\cr
L_2 &=  egja_{1}a_{2} + ehja_{2}b_{1} + fgja_{1}b_{2} + fhjb_{1}b_{2} + cgja_{2} + dgjb_{2} + 2ega_{1} + 2hb_{1}e + 2cg\cr\cr
L_3 &=
	e^{2}gja_{1}a_{2}a_{3} + e^{2}hja_{2}a_{3}b_{1} + efgja_{1}a_{2}b_{3} + efgja_{1}a_{3}b_{2} + efhja_{2}b_{1}b_{3} + efhja_{3}b_{1}b_{2} \cr&\quad + f^{2}gja_{1}b_{2}b_{3} + f^{2}hjb_{1}b_{2}b_{3} + cegja_{1}a_{3} + cegja_{2}a_{3} + cehja_{3}b_{1} + cfgja_{2}b_{3} + degja_{1}b_{3} + degja_{3}b_{2} 
	\cr&\quad + dehjb_{1}b_{3} + dfgjb_{2}b_{3} + c^{2}gja_{3} + cdgjb_{3} + 3e^{2}ga_{1}a_{2} + 3e^{2}ha_{2}b_{1} + 3efga_{1}b_{2} + 3efhb_{1}b_{2} + 3cega_{1} 
	\cr&\quad + 3cega_{2} + 3cehb_{1} + 3degb_{2} + 3c^{2}g
\cr\cr
L_4 &=
e^{3}gja_{1}a_{2}a_{3}a_{4} + e^{3}hja_{2}a_{3}a_{4}b_{1} + e^{2}fgja_{1}a_{2}a_{3}b_{4} + e^{2}fgja_{1}a_{2}a_{4}b_{3} + e^{2}fgja_{1}a_{3}a_{4}b_{2} 
\cr&\quad+ e^{2}fhja_{2}a_{3}b_{1}b_{4} + e^{2}fhja_{2}a_{4}b_{1}b_{3} + e^{2}fhja_{3}a_{4}b_{1}b_{2} + ef^{2}gja_{1}a_{2}b_{3}b_{4} + ef^{2}gja_{1}a_{3}b_{2}b_{4} + ef^{2}gja_{1}a_{4}b_{2}b_{3} 
\cr&\quad+ ef^{2}hja_{2}b_{1}b_{3}b_{4} + ef^{2}hja_{3}b_{1}b_{2}b_{4} + ef^{2}hja_{4}b_{1}b_{2}b_{3} + f^{3}gja_{1}b_{2}b_{3}b_{4} + f^{3}hjb_{1}b_{2}b_{3}b_{4} + ce^{2}gja_{1}a_{2}a_{4} 
\cr&\quad+ ce^{2}gja_{1}a_{3}a_{4} + ce^{2}gja_{2}a_{3}a_{4} + ce^{2}hja_{2}a_{4}b_{1} + ce^{2}hja_{3}a_{4}b_{1} + cefgja_{1}a_{3}b_{4} + cefgja_{1}a_{4}b_{2} + cefgja_{2}a_{3}b_{4}
\cr&\quad + cefgja_{2}a_{4}b_{3} + cefhja_{3}b_{1}b_{4} + cefhja_{4}b_{1}b_{2} + cf^{2}gja_{2}b_{3}b_{4} + de^{2}gja_{1}a_{2}b_{4} + de^{2}gja_{1}a_{4}b_{3} + de^{2}gja_{3}a_{4}b_{2} 
\cr&\quad+ de^{2}hja_{2}b_{1}b_{4} + de^{2}hja_{4}b_{1}b_{3} + defgja_{1}b_{2}b_{4} + defgja_{1}b_{3}b_{4} + defgja_{3}b_{2}b_{4} + defgja_{4}b_{2}b_{3} + defhjb_{1}b_{2}b_{4} 
\cr&\quad+ defhjb_{1}b_{3}b_{4} + df^{2}gjb_{2}b_{3}b_{4} + c^{2}egja_{1}a_{4} + c^{2}egja_{2}a_{4} + c^{2}egja_{3}a_{4} + c^{2}ehja_{4}b_{1} + c^{2}fgja_{3}b_{4} + cdegja_{1}b_{4} 
\cr&\quad+ cdegja_{2}b_{4} + cdegja_{4}b_{2} + cdegja_{4}b_{3} + cdehjb_{1}b_{4} + cdfgjb_{3}b_{4} + d^{2}egjb_{2}b_{4} + 4e^{3}ga_{1}a_{2}a_{3} + 4e^{3}ha_{2}a_{3}b_{1} 
\cr&\quad+ 4e^{2}fga_{1}a_{2}b_{3} + 4e^{2}fga_{1}a_{3}b_{2} + 4e^{2}fha_{2}b_{1}b_{3} + 4e^{2}fha_{3}b_{1}b_{2} + 4ef^{2}ga_{1}b_{2}b_{3} + 4ef^{2}hb_{1}b_{2}b_{3} + c^{3}gja_{4} 
\cr&\quad+ c^{2}dgjb_{4} + 4ce^{2}ga_{1}a_{2} + 4ce^{2}ga_{1}a_{3} + 4ce^{2}ga_{2}a_{3} + 4ce^{2}ha_{2}b_{1} + 4ce^{2}ha_{3}b_{1} + 4cefga_{1}b_{2} + 4cefga_{2}b_{3} 
\cr&\quad+ 4cefhb_{1}b_{2} + 4de^{2}ga_{1}b_{3} + 4de^{2}ga_{3}b_{2} + 4de^{2}hb_{1}b_{3} + 4defgb_{2}b_{3} + 4c^{2}ega_{1} + 4c^{2}ega_{2} + 4c^{2}ega_{3}  
\cr&\quad+ 4c^{2}ehb_{1}+ 4cdegb_{2} + 4cdegb_{3} + 4c^{3}g
\cr
}
$$

\end